\nonstopmode
\documentclass[12pt, a4paper]{amsart}
\usepackage{amscd, amssymb,amsmath,upref}
\usepackage[all]{xy}
\usepackage{amscd}
\usepackage{color}

\allowdisplaybreaks
\begin{document}
\newcommand{\Chi}{\raisebox{2pt}{\ensuremath{\chi}}}
\def\C{\mathbb{C}}
\def\Q{\mathbb{Q}}
\def\N{\mathbb{N}}
\def\Z{\mathbb{Z}}
\def\R{\mathbb{R}}
\def\T{\mathbb{T}}
\def\CS{{\mathcal S}}
\def\CR{{\mathcal R}}
\def\supp{\operatorname{supp}}
\def\Orb{\operatorname{Orb}}
\def\Ind{\operatorname{Ind}}
\def\Aut{\operatorname{Aut}}
\def\End{\operatorname{End}}
\newcommand{\Rep}{\operatorname{Rep}}
\def\Ad{\operatorname{Ad}}
\def\sp{\operatorname{span}}
\def\lsp{\operatorname{span}}
\def\clsp{\overline{\sp}}
\def\To{\mathcal{T}}
\def\Tt{\mathcal{T}}
\def\Ee{\mathcal{E}}
\def\dashind{\operatorname{\!-Ind}}
\def\id{\operatorname{id}}
\def\rt{{\operatorname{rt}}}
\def\lt{\operatorname{lt}}
\def\range{\operatorname{range}}
\def\coker{\operatorname{coker}}
\def\dashind{\operatorname{\!-Ind}}
\def\H{\mathcal{H}}
\def\B{\mathcal{B}}

\def\sign{\operatorname{sign}}

\def\K{\mathcal{K}}
\def\TT{\mathcal{T}}
\def\L{\mathcal{L}}
\def\O{\mathcal{O}}
\renewcommand{\Chi}{\raisebox{2pt}{\ensuremath{\chi}}}
\newcommand{\under}{\backslash}
\newcommand{\qregular}{regular }

\def\tcr{}
\def\aah{}
\def\new{}

\newtheorem{thm}{Theorem}  [section]
\newtheorem{cor}[thm]{Corollary}
\newtheorem{lemma}[thm]{Lemma}
\newtheorem{prop}[thm]{Proposition}
\newtheorem{thm1}{Theorem}
\theoremstyle{definition}
\newtheorem{defn}[thm]{Definition}
\newtheorem{remark}[thm]{Remark}
\newtheorem{example}[thm]{Example}
\newtheorem{examples}[thm]{Examples}
\newtheorem{keyexample}[thm]{Key Example}
\newtheorem{remarks}[thm]{Remarks}
\newtheorem{claim}[thm]{Claim}
\newtheorem{problem}[thm]{Problem}
\newtheorem{conj}[thm]{Conjecture}

\numberwithin{equation}{section}

\title[$C^*$-algebras associated to dilation matrices]
{\boldmath{Purely infinite simple $C^*$-algebras associated to\\ integer dilation matrices}}
\author[Exel]{Ruy Exel}

\address
  {Departamento de Matem\'atica\\
  Universidade Federal de Santa Catarina\\
  88040-900 Florian\'opolis SC\\
  Brazil}
\email{exel@mtm.ufsc.br}

\author[an Huef]{Astrid an Huef}
\address{Department of Mathematics and Statistics\\
University of Otago\\
Dunedin 9054\\
New Zealand}
\email{astrid@maths.otago.ac.nz}

\author[Raeburn]{Iain Raeburn}
\address{School of Mathematics and Applied Statistics\\ The University of Wollongong\\ NSW 2522\\  Australia}
\email{raeburn@uow.edu.au}

\thanks{\new{This research was partially supported by the Australian Research Council and the National Council for Scientific and Technological Development (CNPq) of Brazil.}}

\subjclass[2000]{}

\date{March 5, 2010}
\maketitle

\begin{abstract} \new{Given an $n\times n$ integer matrix $A$ whose eigenvalues
are strictly greater than $1$ in absolute value, let $\sigma_A$ be the
transformation of the $n$-torus $\T^n=\R^n/\Z^n$ defined by $\sigma_A(e^{2\pi
ix})=e^{2\pi iAx}$ for $x\in \R^n$. We study the associated crossed-product
$C^*$-algebra, which is defined using a certain transfer operator for
$\sigma_A$, proving it to be simple and purely infinite and computing its
$K$-theory groups.  } \end{abstract}

\section{Introduction}
Exel has recently introduced a new kind of crossed product for an endomorphism $\alpha$ of a $C^*$-algebra $B$ \cite{E}. The crucial ingredient in his construction is a \emph{transfer operator}, which is a positive linear map $L:B\to B$ satisfying $L(\alpha(a)b)=aL(b)$. In the motivating example, $B=C(X)$, \new{$X$ is a compact Hausdorff space}, $\alpha$ is the endomorphism $\alpha:f\mapsto f\circ\sigma$ associated to a covering map $\sigma:X\to X$, 
and $L$ is defined by
\begin{equation}\label{defL}
L(f)(x)=\frac{1}{|\sigma^{-1}(\{x\})|}\sum_{\sigma(y)=x}f(y).
\end{equation}
Exel's crossed product $B\rtimes_{\alpha,L}\N$ can be constructed in several ways, but here we  view it as the Cuntz-Pimsner algebra $\O(M_L)$ of a right-Hilbert $B$-bimodule $M_L$ constructed from $L$, as discussed in \cite{BR} (see also \S\ref{subsecExel} below). 

We became interested in this circle of ideas when we noticed that the bimodule
$M_L$ associated to the covering map $\sigma:z\mapsto z^N$ of the unit circle
$\T$ plays a key role in work of Packer and Rieffel on projective
multi-resolution analyses \cite{P}--\cite{PR2}. The module elements $m\in M_L$
such that $\langle m,m\rangle$ is the identity of $C(X)$ are precisely the
quadrature mirror filters arising in signal processing and wavelet theory, and
orthonormal bases for $M_L$ are what engineers call ``filter banks with perfect
reconstruction'' (as observed and exploited in \cite{LR2} and \cite{IM}, for
example.) We then noticed further, using results from \cite{EV}, that the
associated crossed product $C(\T)\rtimes_{\alpha_N,L}\N$, 
  where $\alpha_N$ is the endomorphism of
$C(\T^d)$ given by $\sigma$,
  is simple, and
accordingly computed its $K$-theory, finding that $K_0=\Z\oplus(\Z/(N-1)\Z)$ and
$K_1=\Z$. But then we saw this $K$-theory occurring elsewhere, and we gradually
realised that the $C^*$-algebra $C(\T)\rtimes_{\alpha_N,L}\N$ had already been
studied by many authors under other guises. (An almost certainly incomplete list
includes \cite[Example~3]{Deaconu}, \cite[Example~4.1]{KW},
\cite[Appendix~A]{K4} and \cite[Theorem~2.1]{Y}.)

Multiplication by $N$, however, is just one of many dilations of interest in
wavelet theory (see, for example, \cite{S}). Here we consider the covering maps
$\sigma_A$ of $\T^d=\R^d/\Z^d$ induced by integer matrices $A$ whose eigenvalues
$\lambda$ satisfy $|\lambda|>1$, and the crossed products of the associated
systems $(C(\T^d),\alpha_A,L)$,  where $\alpha_A$ is the endomorphism of
$C(\T^d)$ given by $\sigma_A$.

  We show, using results from \cite{EV} and
\cite{K4}, that the crossed products $C(\T^d)\rtimes_{\alpha_A,L}\N$ are simple
and purely infinite, \aah{and} hence by the Kirchberg-Phillips theorem are
classified by their $K$-theory.
  %
The computation of the $K$-theory groups of $C(\T^d)\rtimes_{\alpha_A,L}\N$
therefore has a special significance and one of the main goals of this paper is to
perform precisely this calculation.

Since $C(\T^d)\rtimes_{\alpha_A,L}\N$ is a Cuntz-Pimsner algebra, one should in principle be able to compute its $K$-theory using the exact sequence of \new{{\cite[Theorem~4.8]{Pim}}}, but in practice we were not able to compute some of the homomorphisms in that sequence. So we have argued directly from the six-term exact sequence associated to the Toeplitz algebra of the bimodule $M_L$, and we hope that our computation will be of independent interest. 

Our computation is based on a six-term exact sequence which is valid for any system $(B,\alpha,L)$ for which the bimodule $M_L$ is free as a right Hilbert $B$-module. Using an orthonormal basis for $M_L$, we build a homomorphism $\Omega:B\to M_N(B)$ which has the property that $\Omega\circ \alpha(a)$ is the diagonal matrix $a1_N$ with $N$ copies of $a$ down the diagonal, and which we view as a $K$-theoretic left inverse for $\alpha$. 
When the bimodule is obtained from an integral matrix $A$, as above, this map is closely associated to the classical adjoint of $A$.

We then show that there is an exact sequence
\begin{equation*}
\xymatrix{
K_0(B)\ar[r]^{\id-\Omega_*}&K_0(B)\ar[r]^{j_{B*}\quad}&K_0(\O(M_L))\ar[d]\\
K_1(\O(M_L))\ar[u]&K_1(B)\ar[l]_{\quad \  j_{B*}}&K_1(B)\ar[l]_{\id-\Omega_*}
}
\end{equation*}
in which $j_B$ is the canonical embedding of $B$ in the Cuntz-Pimsner algebra $\O(M_L)$. When $(B,\alpha,L)=(C(\T^d),\alpha_A,L)$, we know from \cite{PR1} that $C(\T^d)_L$ is free, so this exact sequence applies; since we also know from \cite{J} that $K_*(C(\T^d))=K^*(\T^d)$ is isomorphic to the exterior ring generated by a copy of $\Z^d$ in $K^1(\T^d)$, we can in this case compute $\Omega_*$, and derive explicit formulas for $K_i(\O(M_L))$.

\section{Crossed products by endomorphisms}\label{secdefcp}

\subsection{Cuntz-Pimsner algebras} A \emph{right-Hilbert bimodule} over a $C^*$-algebra $B$,
also known as a \emph{correspondence},
is a right Hilbert $B$-module $M$ with a left action of $B$ implemented by a homomorphism $\phi$ of $B$ into the $C^*$-algebra $\L(M)$ of adjointable operators on $M$. In this paper $B$ is always unital, the bimodule $M$ is always \emph{essential} in the sense that $1\cdot m=m$ for $m\in M$, and the bimodule has a \emph{finite Parseval frame} or \emph{quasi-basis}: a finite subset $\{m_j:0\leq j<N\}$ for which we have the \emph{reconstruction formula}
\begin{equation}\label{rp*}
m=\sum_{j =0}^{N-1} m_j\cdot\langle m_j,m\rangle\ \text{ for every $m\in M$.}
\end{equation}
The reconstruction formula implies that 
\begin{equation}\label{leftactfr}
\phi(a)=\sum_{j=0}^{N-1}\Theta_{a\cdot m_j,m_j}\ \text{ for every $a\in B$,}
\end{equation}
and hence that the homomorphism $\phi$ takes values in the algebra $\K(M)$ of compact operators.

The obvious examples of Parseval frames are orthonormal bases:

\begin{lemma}\label{onisframe}
Suppose that $\{m_j:0\leq j<N\}$ are vectors in a right-Hilbert bimodule $M$ over a unital $C^*$-algebra $B$. If the $m_j$ generate $M$ as a Hilbert $B$-module and satisfy $\langle m_j\,,\,m_k\rangle=\delta_{j,k}1_B$, then $\{m_j:0\leq j<N\}$ is a finite Parseval frame for $M$, and $m\mapsto (\langle m_j\,,\,m\rangle)_j$ is an isomorphism of $M$ onto $B^N$.
\end{lemma} 

\begin{proof}
A quick calculation gives the reconstruction formula for $m$ of the form $m_k\cdot b$, and then linearity and continuity give it for arbitrary $m$. For the last assertion, check that $(b_0,\cdots,b_{N-1})\mapsto \sum_jm_j\cdot b_j$ is an inverse.
\end{proof}

\begin{remark}
If $P\in \L(M)$ is a projection and $\{n_j\}$ is an orthonormal basis for $M$, then $\{Pn_j\}$ is a Parseval frame for $P(M)$, and Frank and Larson have shown that every Parseval frame $\{m_j\}$ has this form because $m\mapsto (\langle m_j\,,\,m\rangle)_j$ is an isomorphism of $M$ onto a complemented submodule of $B^N$ \cite[Theorem~5.8]{FR}. However, many interesting bimodules have Parseval frames but are not obviously presented as direct summands of free modules. For example, for a bimodule of the form $C(X)_L$, one can construct a Parseval frame directly using a partition of unity (see, for example, \cite[Proposition~8.2]{EV}). 
\end{remark}

A \emph{Toeplitz representation} of a right-Hilbert bimodule $M$ in a $C^*$-algebra $C$  consists of a linear map $\psi:M\to C$ and a homomorphism $\pi:B\to C$ satisfying $\psi(m)^*\psi(n)=\pi(\langle m\,,\,n\rangle )$ and $\psi(\phi(a)m)=\pi(a)\psi(m)$; we then also have $\psi(m\cdot a)=\psi(m)\pi(a)$.
The \emph{Toeplitz algebra} $\To(M)$ is generated by a universal Toeplitz representation $(i_M,i_B)$ of $M$ (either by theorem \cite{Pim} or by definition \cite{FR}).

The following lemma is implicit in the proof of \cite[Corollary~3.3]{BR}.


\begin{lemma}\label{lem-annoying}
Suppose $M$ is an essential right-Hilbert bimodule over a unital $C^*$-algebra $B$ and $(\psi, \pi)$ is a Toeplitz representation of $M$ on a Hilbert space $\H$. Then the subspace $\pi(1)\H$ is reducing for $(\psi, \pi)$, and
\[
(\psi, \pi)=(\psi_{\pi(1)\H}\oplus 0,\pi_{\pi(1)\H}\oplus 0).
\]
\end{lemma}


\begin{proof}
It is standard that $\pi=\pi_{\pi(1)\H}\oplus 0$, and each $\psi(m)=\psi(1\cdot m)=\pi(1)\psi(m)$ has range in $\pi(1)\H$, so it suffices to show that $h\perp\pi(1)\H$ implies $\psi(m)h=0$. Suppose $h\perp\pi(1)\H$. Then \aah{$\pi(\langle m\,,\, m\rangle)h\in\pi(1)\H$, so that}
\[
\|\psi(m)h\|^2=(\psi(m)h\,|\,\psi(m)h)=(\psi(m)^*\psi(m)h\,|\,h)=(\pi(\langle m\,,\,m\rangle)h\,|\,h)=0.\qedhere
\]
\end{proof}

\begin{remark}\label{convunital}
Lemma~\ref{lem-annoying} implies that the Toeplitz algebra $\To(M)$ is universal for Toeplitz representations $(\psi,\pi)$ in which $\pi$ is unital, and we shall assume from now on that in all Toeplitz representations $(\psi,\pi)$, $\pi$ is unital.
\end{remark}

For every Toeplitz representation $(\psi,\pi)$ of $M$, there is a unique representation $(\psi,\pi)^{(1)}$ of the algebra $\K(M)$ of compact operators on $M$ such that
\[
(\psi,\pi)^{(1)}(\Theta_{m,n})=\psi(m)\psi(n)^*\ \text{ for $m,n\in M$}
\]
(see, for example, \cite[Proposition~1.6]{FR}). When\footnote{As is always the case here; when the left action on the bimodule $M$ contains  non-compact operators, there are several competing definitions of $\O(M)$.} $\phi:B\to \L(M)$ has range in $\K(M)$, we say that $(\psi,\pi)$ is \emph{Cuntz-Pimsner covariant} if $\pi =(\psi,\pi)^{(1)}\circ\phi$, and the \emph{Cuntz-Pimsner algebra} $\O(M)$ is the quotient of $\To(M)$ which is universal for Cuntz-Pimsner covariant representations. The algebra $\O(M)$ is generated by a canonical Cuntz-Pimsner covariant representation $(j_M,j_B)$.

Now we investigate what this all means when $M$ has an orthonormal basis. \new{Compare with \cite[Section~8]{E2} and \cite[Proposition~7.1]{EV} which use quasi-bases.}

\begin{lemma}\label{Cuntzfam}
Suppose that $M$ is an essential right-Hilbert bimodule over a unital $C^*$-algebra $B$, and that $\{m_j:0\leq j<N\}$ is an orthonormal basis for $M$. Let $(\psi,\pi)$ be a Toeplitz representation of $M$. Then:

\smallskip
\textnormal{(1)} $\{\psi(m_j):0\leq j<N\}$ is a Toeplitz-Cuntz family of isometries such that $\sum_{j=0}^{N-1}\psi(m_j)\psi(m_j)^*$ commutes with every $\pi(a)$; and
 
\smallskip
\textnormal{(2)} $(\psi,\pi)$ is Cuntz-Pimsner covariant if and only if $\{\psi(m_j):0\leq j<N\}$ is a Cuntz family.
\end{lemma}

\begin{proof} (1)
The relations $\psi(m_j)^*\psi(m_j)=\pi(\langle m_j\,,\,m_j\rangle)=\pi(1)$ and our convention that $\pi(1)=1$ (see Remark~\ref{convunital}) imply that the $\psi(m_j)$ are isometries. Next,  we fix $a\in B$, let $q:=\sum_{j=0}^{N-1}\psi(m_j)\psi(m_j)^*$, and compute using the reconstruction formula \eqref{rp*}:
\begin{align}\label{compq}
q\pi(a)q&=\sum_{j,k=0}^{N-1}\psi(m_j)\psi(m_j)^*\pi(a)\psi(m_k)\psi(m_k)^*\\
&=\sum_{j,k=0}^{N-1}\psi(m_j)\pi(\langle m_j\,,\, a\cdot m_k\rangle)\psi(m_k)^*\notag\\
&=\sum_{k=0}^{N-1}\Big(\sum_{j=0}^{N-1}\psi(m_j\cdot \langle m_j\,,\, a\cdot m_k\rangle)\Big)\psi(m_k)^*\notag\\
&=\sum_{k=0}^{N-1} \psi( a\cdot m_k)\psi(m_k)^*\notag\\
&=\pi(a)q\notag.
\end{align}
Taking $a=1$ in \eqref{compq} shows that $q^2=q$, and since $q$ is self-adjoint it is a projection. Since each $\psi(m_j)$ is an isometry, each $\psi(m_j)\psi(m_j)^*$ is a projection, and since their sum is a projection, their ranges must be mutually orthogonal. Thus $\{\psi(m_j)\}$ is a Toeplitz-Cuntz family. 
Next we use \eqref{compq} again to see that $q\pi(a)=(\pi(a^*)q)^*=(q\pi(a^*)q)^*=q\pi(a)q=\pi(a)q$, and we have proved (1).

\smallskip\noindent (2)
Suppose that $(\psi,\pi)$ is Cuntz-Pimsner covariant.   Plugging the formula \eqref{leftactfr} for $a=1$ into $(\psi,\pi)^{(1)}(\phi(1))=\pi(1)=1$ shows that $\sum_j\psi(m_j)\psi(m_j)^*=1$, so $\{\psi(m_j)\}$ is a Cuntz family. On the other hand, if $\{\psi(m_j)\}$ is a Cuntz family, then we can deduce from \eqref{leftactfr} that
\[
(\psi,\pi)^{(1)}(\phi(a))=\sum_{j=0}^{N-1}\psi(m_j)\psi(a^*\cdot m_j)^*=\sum_{j=0}^{N-1}\psi(m_j)\psi(m_j)^*\pi(a)=\pi(a),
\]
and $(\psi,\pi)$ is Cuntz-Pimsner covariant.
\end{proof}

\subsection{Exel systems and crossed products}\label{subsecExel}

Let $\alpha$ be an endomorphism of a unital $C^*$-algebra $B$.   A \emph{transfer operator} $L$ for $(B,\alpha)$ is a positive linear map $L:B\to B$ such that $L(\alpha(a)b)=aL(b)$ for all $a,b\in B$. We call the triple $(B,\alpha, L)$ an \emph{Exel system}. 

Given an Exel system $(B,\alpha, L)$, we construct a right-Hilbert $B$-module $M_L$ over $B$ as in \cite{E} and \cite{BR}.  Let $B_L$ be a copy of the underlying vector space of $B$.  Define a right action of $a\in B$ on $m\in B_L$ by $m\cdot a=m\alpha(a)$, and a $B$-valued pairing on $B_L$ by
\[
\langle m\,,\,n\rangle=L(m^*n)\quad\text{for $m,n\in B_L$}.
\]
Modding out by $\{m:\langle m\,,\,m\rangle=0\}$ and completing yields a right Hilbert $B$-module $M_L$. The action of $B$ by left multiplication  on $B_L$ extends to an action of $B$ by adjointable operators on $M_L$ which is implemented by a unital homomorphism $\phi:B\to \L(M_L)$, and and thus makes $M_L$ into a right-Hilbert bimodule over $B$.

Exel's crossed product is constructed in two stages. First he forms a Toeplitz algebra $\To(B,\alpha,L)$, which is isomorphic to $\To(M_L)$ (see \cite[Corollary~3.2]{BR}). Then the crossed product $B\rtimes_{\alpha,L}\N$ is the quotient of $\To(M_L)$ by the ideal generated by the elements
\[
i_B(a)-(i_{M_L},i_B)^{(1)}(\phi(a))\ \text{ for $a\in K_\alpha:=\phi^{-1}(\K(M_L))
\cap\overline{B\alpha(B)B}$}
\]
(see \cite[Lemma~3.7]{BR}). When $M_L$ has a finite Parseval frame and the projection $\alpha(1)$ is full, we have $\phi^{-1}(\K(M_L))=B=K_\alpha$, and $B\rtimes_{\alpha,L}\N$ is the Cuntz-Pimsner algebra $\O(M_L)$.

For us, the main examples of Exel systems come from surjective endomorphisms $\sigma$ of a compact group $K$ with finite kernel: the corresponding Exel system $(C(K),\alpha ,L)$ has $\alpha(f)=f\circ\sigma$ and $L$ defined by averaging over the fibres of $\sigma$, as in \eqref{defL}. The next lemma is a mild generalisation of \cite[Proposition~1]{PR1}.

\begin{lemma}\label{ex-compactgroup}
Suppose that $\sigma:K\to K$ is a surjective endomorphism of a compact abelian group $K$ with $N:=|\ker\sigma|<\infty$, and $(C(K),\alpha,L)$ is the corresponding Exel system. Then the norm on $C(K)_L$ defined by the inner product is equivalent to the usual sup-norm, and $C(K)_L$ is complete. It has an orthonormal basis $\{m_j:0\leq j<N\}$. 
\end{lemma}

\begin{proof}
The assertions about the norm and the completeness are proved in \cite[Lemma~3.3]{LR2}, for example.  Since $\gamma\mapsto \gamma|_{\ker\sigma}$ is surjective and $|(\ker\sigma)^\wedge|=|\ker\sigma|=N$, we can find a subset $\{\gamma_i:0\leq i<N\}$ of $\widehat K$ such that $\{\gamma_i|_{\ker\sigma}\}$ is all of $(\ker\sigma)^\wedge$. Then
\begin{align*}
\langle \gamma_i,\gamma_j\rangle_L(k)&=\frac{1}{N}\sum_{\sigma(l)=k}\overline{\gamma_i(l)}\gamma_j(l)\\
&=\frac{1}{N}\sum_{\zeta\in\ker \sigma}\overline{\gamma_i(\zeta l_0)}\gamma_j(\zeta l_0)\text{ for any fixed $l_0$ such that $\sigma(l_0)=k$}\\
&=\frac{1}{N}\overline{\gamma_i(l_0)}\gamma_j(l_0)\sum_{\zeta\in\ker \sigma}\big(\overline{\gamma_i}\gamma_j\big)(\zeta).
\end{align*}
If $i\not= j$, then $(\gamma_i^{-1}\gamma_j)|_{\ker\sigma}$ is a nontrivial character of $\ker \sigma$, and its range is a nontrivial subgroup of $\T$, so the sum vanishes. If $i=j$, then the sum is $N$. So $\{\gamma_j\}$ is orthonormal. 

We still need to see that $\{\gamma_j\}$ generates $C(K)_L$ as a Hilbert module. The Stone-Weierstrass theorem implies that the characters of $K$ span a dense $*$-subalgebra of $C(K)$, and hence by the equivalence of the norms, they also span a dense subspace of $C(K)_L$. So it suffices to show that each $\gamma\in \widehat K$ is in the submodule generated by $\{\gamma_j\}$. Since $(\ker\sigma)^\wedge=\{\gamma_j\}$, there exists $j$ such that $\gamma|_{\ker\sigma}=\gamma_j$. Then $\gamma_j^{-1}\gamma$ vanishes on $\ker \sigma$, and there is a character $\chi$ such that  $\gamma_j^{-1}\gamma=\chi\circ\sigma$. This equation unravels as \aah{$\gamma=\gamma_j(\chi\circ\sigma)=\gamma_j\alpha(\chi)=\gamma_j\cdot\chi$}, so it implies that $\gamma$ belongs to the submodule generated by $\{\gamma_j\}$. 
\end{proof}

\begin{example}
Suppose that $A\in M_d(\Z)$ is an integer matrix with $|\det A|>1$, and $\sigma_A$ is the endomorphism of $\T^d$ given by $\sigma_A(e^{2\pi ix})=e^{2\pi iAx}$ for $x\in \R^d$. Then $\sigma_A$ is surjective (because $A:\R^d\to \R^d$ is), and $\ker \sigma_A$ has $N:=|\det A|$ elements. A function $m\in C(\T^d)_L$ such that $\langle m\,,\,m\rangle(z)=1$ for all $z$ is called a \emph{quadrature mirror filter} for dilation by $A$, and an orthonormal basis  for $C(\T^d)_L$ is a \emph{filter bank}. Lemma~\ref{ex-compactgroup} says that for every $A$, filter banks exist. 
\end{example}

\begin{remark} 
Although the dilation matrices $A$ are of great relevance to wavelets, the filter banks we constructed in the proof of Lemma~\ref{ex-compactgroup} are not the kind which are useful for the construction of wavelets. There one wants the first filter $m_0$ to be low-pass, which means roughly that $m_0(1)=N^{1/2}$, $m_0$ is smooth near $1$, and $m_0$ does not vanish on a sufficiently large neighbourhood of $1$; for the basis in Lemma~\ref{ex-compactgroup}, we have $|m_0(z)|=1$ for all $z$, and $m_0$ is \emph{all-pass}. The \emph{matrix completion problem} considered in \cite{PR1} asks whether, given a low-pass filter $m_0$, one can find a filter bank $\{m_j\}$ which includes the given $m_0$. This amounts to asking that the submodule $m_0^\perp:=\{m\in C(\T^d)_L:\langle m\,,\,m_0\rangle=0\}$ is free. In  \cite[\S4]{PR1}, Packer and Rieffel show by example that it need not be free if $|\det A|>2$ and $d>4$. Of course, since $m_0^\perp$ is a direct summand of a free module, it always has a Parseval frame.
\end{remark}

When $\alpha$ is the endomorphism of $C(K)$ coming from a surjective endomorphism $\sigma$ of $K$, we know from Lemma~\ref{ex-compactgroup} that $M_L=C(K)_L$ admits an orthonormal basis, and the associated endomorphism $\alpha:f\mapsto f\circ\sigma$ is unital, so $\alpha(1)=1$ is certainly full. Thus for the systems of interest to us, Exel's crosed product $B\rtimes_{\alpha,L}\N$ is isomorphic to the Cuntz-Pimsner algebra $\O(M_L)$. We will use this identification without comment.

\section{The six-term exact sequence}

We assume throughout this section that $(B,\alpha,L)$ is an Exel system and that $\{m_j:0\leq j<N\}$ is a Parseval frame for $M_L$. We write $Q$ for the quotient map from $\To(M_L)\to\O(M_L)$, and $(\psi,\pi)$ for the universal Toeplitz covariant representation of $M_L$ in $\To(M_L)$. 

To construct our exact sequence for $K_*(\O(M))$, we analyse the six-term exact sequence
\begin{equation}\label{eq-les2}
\xymatrix{
K_0(\ker Q)\ar[r]^{\iota_*}&K_0(\To(M_L))\ar[r]^{Q_*}&K_0(\O(M_L))\ar[d]_{\delta_0}\\
K_1(\O(M_L))\ar[u]_{\delta_1}&K_1(\To(M_L))\ar[l]_{Q_*}&K_1(\ker Q).\ar[l]_{\iota_*}
}
\end{equation}
We begin by recalling from \cite[Theorem~4.4]{Pim} that the homomorphism $\pi:B\to \To(M_L)$ induces an isomorphism of $K_i(B)$ onto $K_i(\To(M_L))$, so we can replace $K_i(\To(M_L))$ with $K_i(B)$ provided we can identify the maps. Next we introduce our ``$K$-theoretic left inverse'' for $\alpha$, and then we will work towards showing that $B$ is a full corner in $\ker Q$, so that we can replace $K_i(\ker Q)$ with $K_i(B)$.
 
Part (3) of the next result will not be used in this section; it is included here because it shows how $\Omega$ relates to $\alpha$, and gives a hint of why we view it as a ``$K$-theoretic left inverse'' for $\alpha$.

\begin{lemma}\label{lem-Omega} Define $\Omega: B\to M_N(B)$ by $\Omega(a)=(\langle m_j\,,\, a\cdot m_k\rangle)_{j,k}$.  Then

\smallskip
\textnormal{(1)} $\Omega$ is a homomorphism of $C^*$-algebras;

\smallskip
\textnormal{(2)} $\Omega$ is unital if and only if $\{m_j:0\leq j<N\}$ is an orthonormal basis;

\smallskip
\textnormal{(3)} if $B$ is commutative and $\{m_j:0\leq j<N\}$ is orthonormal, then $\Omega(\alpha(a))$ is the diagonal matrix $a1_N$ with  \aah{diagonal entries $a$.}
\end{lemma}

\begin{proof}
For (1), we let $a,b\in B$ and compute: first
\begin{align*}
(\Omega(a)\Omega(b))_{j,k}
&=\sum_{l=0}^{N-1}\langle m_j\,,\, a\cdot m_l\rangle \langle m_l\,,\, b\cdot m_k\rangle\\
&=\Big\langle m_j\,,\, a\cdot \Big(\sum_{l=0}^{N-1}m_l\cdot \langle m_l\,,\, b\cdot m_k\rangle\Big)\Big\rangle \\
&=\langle m_j\,,\, a\cdot (b\cdot m_k)\rangle\\
&=\Omega(ab)_{j,k},
\end{align*}
and then
\[
\Omega(a^*)=(\langle m_j\,,\, a^*\cdot m_k\rangle)_{j,k}=(\langle  a\cdot m_j\,,\, m_k\rangle)_{j,k}=(\langle m_k\,,\, a\cdot m_j\rangle^*)_{j,k}=\Omega(a)^*.
\]

Part (2) is easy. For (3), we let $q_L:B_L\to M_L$ be the quotient map, and consider $m=q(b)\in q(B_L)$. Then commutativity of $B$ gives
\[m\cdot a=q(b\cdot a)=q(b\alpha(a))=q(\alpha(a)b)=\alpha(a)\cdot q(b)=\alpha(a)\cdot m,
\]
and this formula extends to $m\in M_L$ by continuity. Thus
\begin{align*}
\Omega(\alpha(a))_{j,k}
=\langle m_j\,,\,\alpha(a)\cdot m_k\rangle=\langle m_j\,,\,m_k\cdot a\rangle
=\langle m_j\,,\,m_k\rangle a=\delta_{j,k}a,
\end{align*}
as required.
\end{proof}

To describe $\ker Q$, we need some standard notation. We write $M_L^{\otimes i}$ for the $i$-fold internal tensor product $M_L\otimes_B \cdots\otimes_B M_L$, which is itself a right-Hilbert bimodule over $B$. There is a Toeplitz representation $(\psi^{\otimes i},\pi)$ of $M_L^{\otimes k}$ in $\To(M_L)$ such that $\psi^{\otimes k}(\xi)=\prod_{i=1}^i\psi(\xi_i)$ for elementary tensors $\xi=\xi_1\otimes\cdots\otimes\xi_k$ in $M_L^{\otimes k}$ (see \cite[Proposition~1.8]{FR}, for example). By convention, we set $M_L^{\otimes 0}:=B$ and $\psi^{\otimes 0}:=\pi$. Then from \cite[Lemma~2.4]{FR} we have
\begin{equation}\label{spanTo}
\To(M_L)=\clsp\{\psi^{\otimes k}(\xi)\psi^{\otimes l}(\eta)^*:k,l\geq 0, \xi\in M_L^{\otimes k}, \eta\in M_L^{\otimes l}\}.
\end{equation}
We also recall from Lemma~\ref{Cuntzfam}(1) that the element $q:=\sum_{j=0}^{N-1}\psi(m_j)\psi(m_j)^*$ of $\To(M_L)$ is a projection which commutes with every $\pi(a)$.

\begin{lemma}\label{lem-kerQ}
With the preceding notation, we have
\begin{enumerate}
\item  $1-q=1-\sum_{j=0}^{N-1}\psi(m_j)\psi(m_j)^*$ is a full projection in $\ker Q$;
\smallskip
\item $(1-q)\psi^{\otimes k}(\xi)=0$ for all $\xi\in M_L^{\otimes k}$ with $k\geq 1$; and
\smallskip
\item
$\ker Q=\clsp\{\psi^{\otimes k}(\xi)(1-q)\psi^{\otimes l}(\eta)^*:k,l\geq 0, \xi\in M_L^{\otimes k}, \eta\in M_L^{\otimes l}\}$.
\end{enumerate}
\end{lemma}

\begin{proof}
(1) The reconstruction formula implies that
$\phi(a)=\sum_{j=0}^{N-1}\Theta_{a\cdot m_j,m_j}$, and so 
\begin{equation}\label{formfor^(1)}
(\psi,\pi)^{(1)}(\phi(a))=\sum_{j=1}^{N-1}\psi(a\cdot m_j)\psi(m_j)^*=\pi(a)q.
\end{equation}
This implies in particular that
\[
Q(1-q)=Q(\pi(1)-\pi(1)q)=Q\big(\pi(1)-(\psi,\pi)^{(1)}(\phi(1))\big)=0, 
\]
so $1-q$ belongs to $\ker Q$. Since $\ker Q$ is by definition the ideal in $\To(M_L)$ generated by the elements $\pi(a)-(\psi,\pi)^{(1)}(\phi(a))$ for $a\in B$, \eqref{formfor^(1)} also implies that $\ker Q$ is generated by the elements $\pi(a)(1-q)$, and hence by the single element $1-q$. This says precisely that the projection $1-q$ is full.

(2) First we consider $m\in M_L^{\otimes 1}=M_L$. The reconstruction formula gives
\[
q\psi(m)
=\sum_{j=0}^{N-1}\psi(m_j)\psi(m_j)^*\psi(m)
=\sum_{j=0}^{N-1}\psi(m_j\cdot \langle m_j\,,\, m\rangle )=\psi(m),
\]
so $(1-q)\psi(m)=0$. Now for $k>1$ and for an elementary tensor $\xi=\xi_1\otimes\cdots \otimes\xi_k$, we have $(1-q)\psi^{\otimes k}(\xi)=(1-q)\big(\prod_{i=1}^i\psi(\xi_i)\big)=0$, and the result extends to arbitrary $\xi \in M^{\otimes k}$ by linearity and continuity.

(3) In view of part (2), we can deduce from \eqref{spanTo} that  $\ker Q=\To(M_L)(1-q)\To(M_L)$ is spanned by the elements of the form
\[\psi^{\otimes k}(\xi)\pi(a)^*(1-q)\pi(b)\psi^{\otimes l}(\eta)^*=
\psi^{\otimes k}(\xi\cdot a)(1-q)\psi^{\otimes l}(\eta\cdot b^*)^*\]
for $\xi\in M^{\otimes k}$, $\eta\in M^{\otimes l}$ and $a,b\in B$, which gives (3).
\end{proof}

\begin{lemma}\label{lem-rho}
There is a homomorphism $\rho:B\to\ker Q$ such that $\rho(a)=\pi(a)(1-q)$, and  $\rho$ is an isomomorphism of $B$ onto $(1-q)\ker Q(1-q)$.
\end{lemma}

\begin{proof} Lemma~\ref{lem-kerQ} says that $\pi(a)(1-q)$ belongs to $\ker Q$ and Lemma~\ref{lem-Omega} says that $q$ commutes with every $\pi(a)$, so there is a homomorphism $\rho:B\to (1-q)\ker Q(1-q)\subset\ker Q$ such that $\rho(a)=\pi(a)(1-q)$.  From parts (2) and (3) of Lemma~\ref{lem-kerQ} we get:
\begin{align*}
(1-q)\ker Q(1-q)&=\clsp\{ (1-q)\psi^{\otimes k}(\xi)(1-q)\psi^{\otimes j}(\eta)^*(1-q): k,l\geq 0 \}\\
&=\clsp\{ (1-q)\pi(a)(1-q)\pi(b)(1-q): a,b\in B \}\\
&=\clsp\{ (1-q)\pi(ab): a,b\in B \},
\end{align*}
which is precisely the range of $\rho$. So $\rho$ is surjective.

To see that $\rho$ is injective we choose a faithful representation $\pi_0:B\to B(\H)$ and consider the Fock representation $(\psi_F,\pi_F)$ of $M_L$ induced from $\pi_0$, as described in \cite[Example~1.4]{FR}. The underlying space of this Fock representation is $F(M_L)\otimes_B\H:=\bigoplus_{k\geq 0}(M_L^{\otimes k}\otimes_B\H)$; $B$ acts diagonally on the left, and $M_L$ acts by creation operators. The crucial point for us is that each $\psi_F(m)^*$ is an annihilation operator which vanishes on the subspace $B\otimes_B\H=M_L^{\otimes 0}\otimes_B\H$ of $F(M_L)\otimes_B\H$. 

Now suppose that $a\in B$. Then
\[
0=\psi_F\times\pi_F(\rho(a))= \psi_F\times\pi_F(\pi(a)(1-q))=\pi_F(a)\Big(1-\sum_{j=0}^{N-1}\psi_F(m_j)\psi_F(m_j)^*\Big).
\]
Since $\psi(m_j)^*$ vanishes on $B\otimes_B \H$, we have
\begin{align*}
\rho(a)=0&\Longrightarrow \pi_F(a)\Big(1-\sum_{j=0}^{N-1}\psi_F(m_j)\psi_F(m_j)^*\Big)(1\otimes_B h)=0\ \text{ for all $h\in \H$}\\
&\Longrightarrow \pi_F(a)(1\otimes_B h)=0\ \text{ for all $h\in \H$}\\
&\Longrightarrow a\otimes_B h=0\ \text{ for all $h\in \H$}\\
&\Longrightarrow \pi_0(a)h=0\ \text{ for all $h\in \H$,}
\end{align*}
which implies that $a=0$ because $\pi_0$ is faithful.
\end{proof}

Lemma~\ref{lem-rho} implies we can replace $K_i(\ker Q)$ in \eqref{eq-les2} by $K_i(B)$, as claimed. Now we need to check what this replacement does to the map $\iota_*$.

\begin{prop}\label{prop-square}
The following diagram commutes for $i=0$:
\begin{equation}\label{diagram**}
\xymatrix{
K_i(B)\ar[r]^{\id-\Omega_*}\ar[d]_{\rho_*} &K_i(B)\ar[d]_{\pi_*}\\
K_i(\ker Q)\ar[r]_{\iota_*}&K_i(\To(M_L))
}
\end{equation}
If  $\{m_j:0\leq j<N\}$ is orthonormal then the diagram also commutes for $i=1$.
\end{prop}

Since $\Omega:B\to M_N(B)$, the $\Omega_*$ in the diagram is really the composition of $\Omega_*:K_i(B)\to K_i(M_N(B))$ with the isomorphism $K_i(M_N(B))\to K_i(B)$; the latter is induced by the map which views an element in $M_r(M_N(B))$ as an element of $M_{rN}(B)$. 

The proof needs two standard lemmas. The first says, loosely, that if we rewrite an $r\times r$ matrix of $N\times N$ blocks as an $N\times N$ matrix of $r\times r$ blocks, then the resulting $rN\times rN$ matrices are unitarily equivalent. We agree that this can't be a surprise to anyone, and we apologise for failing to come up with more elegant notation.

\begin{lemma}\label{tocome}
Suppose that $B$ is a $C^*$-algebra, $r\geq 1$ and $N\geq 2$ are integers, and
\[
\{b_{j,s;k,t}:0\leq j,k<N\text{ and }0\leq s,t<r\}
\]
is a subset of $B$. For $m,n$ satisfying $0\leq m,n<rN-1$, we define 
\begin{align*}
c_{m,n}&=b_{j,s;k,t}\text{ where $m=sN+j$ and $n=tN+k$, and}\\
d_{m,n}&=b_{j,s;k,t}\text{ where $m=jr+s$ and $n=kr+t$.}
\end{align*}
Then there is a scalar unitary permutation matrix $U$ such that the matrices $C:=(c_{m,n})$ and $D:=(d_{m,n})$ are related by $C=UDU^*$.
\end{lemma}

\begin{proof}
For $0\leq p,q<rN-1$, we define
\[
u_{p,q}=\begin{cases}
1&\text{if there exist $k$, $t$ such that $p=tN+k$ and $q=kr+t$}\\
0&\text{otherwise.}
\end{cases}
\]
Each row and column contain exactly one $1$, so $U:=(u_{p,q})$ is a scalar permutation matrix, and we can verify that both $(CU)_{m,q}$ and $(UD)_{m,q}$ are equal to $b_{j,s;k,t}$ where $m=sN+j$ and $q=kr+t$, so $CU=UD$.
\end{proof}

\begin{lemma}\label{lem-isometrytrick}
Suppose that $S$ is an isometry in a unital $C^*$-algebra $B$. Then
\[
U:=\begin{pmatrix}
S&1-SS^*\\
0&S^*
\end{pmatrix}
\]
is a unitary element of $M_2(B)$
 and its class in $K_1(B)$ is the identity.
\end{lemma}
\begin{proof} A straightforward calculation shows that $U$ is unitary.

Let $\To=C^*(v)$ be the Toeplitz algebra.
By Coburn's Theorem \cite{Coburn} there is  a homomorphism $\pi_S:\To\to B$ such that $\pi_S(v)=S$.  Since $K_1(\To)=0$ (see, for example, \cite[Remark~11.2.2]{WO}),
\[
\left[\begin{pmatrix}
v&1-vv^*\\
0&v^*
\end{pmatrix}\right]=[1] \text{\ in $K_i(\To)$,}
\]
and hence
\[
\left[\begin{pmatrix}
S&1-SS^*\\
0&S^*
\end{pmatrix}\right]
=(\pi_S)_*\left( \left[\begin{pmatrix}
v&1-vv^*\\
0&v^*
\end{pmatrix}\right] \right)=(\pi_S)_*([1])=[1] \text{\ in $K_i(B)$.}\qedhere
\]
\end{proof}

\begin{proof}[Proof of Proposition~\ref{prop-square}]
  We start with $i=0$. Let $a=(a_{s,t})$ be a projection in $M_r(B)$. For $\pi:A\to B$, we  write $\pi_r$ for the induced homomorphism of $M_r(A)$ into $M_r(B)$. Then we have
\begin{align*}
\rho_*([a])&=[(\rho(a_{r,s})]=\big[(\pi(a_{s,t})(1-q))\big]=\big[(\pi(a_{s,t}))(1-q)1_r)\big]\\&=[\pi_r(a)((1-q)1_r)]=[\pi_r(a)]-[(\pi_r(a)(q1_r))],\quad\text{and}\\
\pi_*\circ(\id-\Omega_*)([a])&=[\pi_r(a)]-\pi_{*}\circ\Omega_*([a]),
\end{align*}
so it suffices to show that $[\pi_r(a)(q1_r)]=\pi_{*}\circ\Omega_*([a])$ in $K_0(\To(M_L))$. The class $\pi_{*}\circ\Omega_*([a])$ appears as the class of the $r\times r$ block matrix $\pi_{rN}(\Omega_r(a))$ whose $(s,t)$ entry is the $N\times N$ block $\big(\pi(\langle m_j\,,\, a_{s,t}\cdot m_k\rangle)\big)_{j,k}$. In other words, with $b_{j,s;k,t}=\pi(\langle m_j\,,\, a_{s,t}\cdot m_k\rangle)$, \aah{the matrix} $\pi_{rN}(\Omega_r(a))$ is the matrix $C=(c_{m,n})$ in Lemma~\ref{tocome}.

We now consider the matrix $T$ in $M_N(M_r(\To(M_L)))$ defined by 
\begin{equation}
\label{eq-T}T=\begin{pmatrix}
\psi(m_0)1_r&\cdots&\psi(m_{N-1})1_r\\
0_r&\cdots&0_r\\
\vdots&\cdots&\vdots
 \end{pmatrix}. 
 \end{equation}
Computations show that $TT^*=(q1_r)\oplus  0_{r(N-1)}$, and since $\pi_r(a)$ is a projection which commutes with $q1_r$, we deduce that $(\pi_r(a)\oplus  0_{r(N-1)})T$ is a partial isometry which implements a Murray-von Neumann equivalence between $T^*(\pi_r(a)\oplus  0_{r(N-1)})T$ and $\big(\pi_r(a)\oplus 0_{r(N-1)}\big)TT^*=(\pi_r(a)(q1_r))\oplus 0_{r(N-1)}$. Thus we have
\[
[\pi_r(a)(q1_r)]=\big[\pi_r(a)(q1_r)\oplus 0_{r(N-1)}\big]=\big[T^*(\pi_r(a)\oplus  0_{r(N-1)})T\big].
\]
Another computation shows that the $(j,k)$ entry of $T^*(\pi_r(a)\oplus  0_{r(N-1)})T$ is the $r\times r$ matrix $\big(\pi_r(\langle m_j\,,\, a_{s,t}\cdot m_k\rangle1_r)\big)_{s,t}$.  Thus with the same choice of $b_{j,s;k,t}=\pi(\langle m_j\,,\, a_{s,t}\cdot m_k\rangle)$, $T^*(\pi_r(a)\oplus  0_{r(N-1)})T$ is the matrix $D=(d_{m,n})$ in Lemma~\ref{tocome}. Since unitarily equivalent projections have the same class in $K_0$, we can therefore deduce from Lemma~\ref{tocome} that
\begin{equation}\label{changeshape}
[\pi_r(a)(q1_r)]=\big[T^*(\pi_r(a)\oplus  0_{r(N-1)})T\big]=\big[\pi_{rN}(\Omega_r(a))\big]=\pi_*\circ\Omega_*([a])].
\end{equation}
Thus Diagram~\ref{diagram**} commutes when $i=0$.

Now consider $i=1$, where we assume in addition that $\{m_j\}$ is orthonormal.  Let $u$ be a unitary in $M_r(B)$. To compute $\rho_*:K_1(B)\to K_1(\ker Q)$ we observe that $\rho$ is the composition  of a unital isomorphism of $B$ onto $(1-q)\ker Q (1-q)$, which takes $[u]$ to $[\rho_r(u)]=[\pi_r(u)((1-q)1_r)]$, with the inclusion of $(1-q)\ker Q (1-q)$ as a full corner in the non-unital algebra $\ker Q$, which takes 
$[\pi_r(u)((1-q)1_r)]$ to $[\pi_r(u)((1-q)1_r)+q1_r]\in K_1((\ker Q)^+)=K_1(\ker Q)$.
On the other hand,
\[
\pi_*\circ(\id-\Omega_*)([u])=[\pi_r(u)]-[\pi_{rN}\circ\Omega_r(u)].
\]
So we need to show that
\begin{equation}\label{eq-strategy}
[(\pi_r(u)((1-q)1_r)+ q1_r)\oplus 1_{r(N-1)}]=[\pi_r(u)\oplus 1_{r(N-1)}]-[\pi_{rN}\circ\Omega_r(u)]
\end{equation}
in $K_1(\To(M_L))$. To this end, we note that the left-hand side of \eqref{eq-strategy} is unchanged by pre- or post-multiplying by any invertible matrix $C\in M_{2rN}(\To(M_L))$ whose $K_1$ class is $1$. In particular, we can do this when $C$ is:
\begin{itemize}
\item a  unitary of the form
\[
C=\begin{pmatrix}
S&1-SS^*\\
0&S^*
\end{pmatrix}\] where $S\in M_{rN}(\To(M_L))$ is an isometry (see Lemma~\ref{lem-isometrytrick});
\item
an upper- or lower-triangular matrix of the form
\[
C=\begin{pmatrix}
1&A\\
0&1
\end{pmatrix}\quad\text{or}\quad C=\begin{pmatrix}
1&0\\
A&1
\end{pmatrix}
\]
(which are connected to $1_{2rN}$ via $t\mapsto\left( \begin{smallmatrix}
1&tA\\
0&1
\end{smallmatrix}\right)$ and its transpose);
\item
 any constant invertible matrix $C$ in $M_{2rN}(\C)$ (because $GL_{2rN}(\C)$ is connected); this implies that we can perform row and column operations without changing the class in $K_1$.
 \end{itemize}

Since $\{m_j\}$ is an orthonormal basis, the matrix $T$ defined at \eqref{eq-T} is an isometry in $M_{rN}(\To(M_L))$. Thus
\begin{align*}
\big[\big(\pi_r(u)&((1-q)1_r)+ q1_r\big)\oplus 1_{r(N-1)}\big]\\
&=\left[
\begin{pmatrix}(\pi_r(u)((1-q)1_r)+ q1_r)\oplus 1_{r(N-1)}&0_{rN} \\0_{rN}& 1_{rN} 
\end{pmatrix}\right]
\left[\begin{pmatrix}
T&1_{rN}-TT^*\\
0_{rN}&T^*
\end{pmatrix}\right]\\
&=\left[\begin{pmatrix}(\pi_r(u)((1-q)1_r)+ q1_r)\oplus 1_{r(N-1)}&0_{rN} \\0_{rN}& 1_{rN} \end{pmatrix}\right]
 \left[\begin{pmatrix}
T&(1-q)1_r\oplus 1_{r(N-1)} \\
0_{rN}&T^*
\end{pmatrix}\right]\\
&= \left[
\begin{pmatrix}
\aah{\big((\pi_r(u)((1-q)1_r)+ q1_r)\oplus 1_{r(N-1)}\big)
T}
&
\pi_r(u)((1-q)1_r)\oplus 1_{r(N-1)}
\\0_{rN}&T^*
\end{pmatrix}
\right],\\
\intertext{\aah{which, since $(1-q)\psi(m_i)=0$ by Lemma~\ref{lem-kerQ}(2), is }}
&=\left[\begin{pmatrix}
T &
\pi_r(u)((1-q)1_r)\oplus 1_{r(N-1)}
\\
0_{rN} &T^*
\end{pmatrix}\right]\\
&=\left[\begin{pmatrix}
T &
\pi_r(u)((1-q)1_r)\oplus 1_{r(N-1)}\\
0_{rN}&T^*
\end{pmatrix}\right]\left[ \begin{pmatrix}1_{rN}&T^*\big(\pi_r(u)\oplus 1_{r(N-1)}  \big) \\0_{rN}&1_{rN}\end{pmatrix} \right]\\
&=\left[\begin{pmatrix}T&\pi_r(u)\oplus 1_{r(N-1)}\\0_{rN}& T^* \end{pmatrix} \right]\\
\intertext{since $TT^*=q1_r\oplus 0_{r(N-1)}$ and $(q1_r)\pi_r(u)=\pi_r(u)(q1_r)$. 
By an elementary row operation this is}
&=\left[\begin{pmatrix}\pi_r(u)\oplus 1_{r(N-1)}&T\\ T^*&0_{rN} \end{pmatrix} \right]\\
&=\left[\begin{pmatrix}\pi_r(u)\oplus 1_{r(N-1)}&T\\ T^*& 0_{rN}\end{pmatrix} \right]\left[\begin{pmatrix}1_{rN}&-\big(\pi_r(u^{-1})\oplus 1_{r(N-1)}\big)T\\ 0_{rN}& 1_{rN}\end{pmatrix} \right]\\
&=\left[\begin{pmatrix}
\pi_r(u)\oplus 1_{r(N-1)}&0_{rN}\\T^*&-T^*\big(\pi_r(u^{-1})\oplus 1_{r(N-1)}  \big)T
\end{pmatrix} \right]\\
&=\left[
\begin{pmatrix}1_{rN}&0_{rN}\\-T^*\big(\pi_r(u^{-1})\oplus 1_{r(N-1)}\big) & 1_{rN}\end{pmatrix}
\right]
\left[\begin{pmatrix}
\pi_r(u)\oplus 1&0_{rN}\\T^*&-T^*\big(\pi_r(u^{-1})\oplus 1  \big)T
\end{pmatrix} \right]\\
&=\left[\begin{pmatrix}
\pi_r(u)\oplus 1_{r(N-1)}&0_{rN}\\0_{rN}&-T^*\big(\pi_r(u^{-1})\oplus 1_{r(N-1)}  \big)T
\end{pmatrix} \right]\left[\begin{pmatrix}1_{rN}&0_{rN}\\0_{rN}&-1_{rN} \end{pmatrix}\right]\\
&=[\pi_r(u)\oplus 1_{r(N-1)}]+[T^*\big(\pi_r(u^{-1})\oplus 1_{r(N-1)}  \big)T].
\end{align*}
Now we recall from the argument in the second paragraph (see \eqref{changeshape}) that
\[
\big[T^*\big(\pi_r(u^{-1})\oplus 1_{r(N-1)}  \big)T\big]=\big[\pi_{rN}(\Omega_r(u^{-1}))\big]=-[\pi_{rN}\circ\Omega_r(u)],
\]
and we see that we have proved what we wanted.
\end{proof}

\begin{thm}\label{thm-stes}
Let $(B,\alpha,L)$ be an Exel system with $B$ unital and separable, and suppose that $M_L$ has an orthonormal basis $\{m_j\}_{j=0}^{N-1}$.  Let $(j_{M_L},j_B)$ be the canonical Cuntz-Pimsner covariant representation of $M_L$ in $\O(M_L)$. Then there is an exact sequence
\begin{equation}\label{eq-les}
\xymatrix{
K_0(B)\ar[r]^{\id-\Omega_*}&K_0(B)\ar[r]^{j_{B*}\quad}&K_0(\O(M_L))\ar[d]^{\rho_*^{-1}\circ\delta_0}\\
K_1(\O(M_L))\ar[u]^{\rho_*^{-1}\circ\delta_1}&K_1(B)\ar[l]_{\quad \  j_{B*}}&K_1(B).\ar[l]_{\id-\Omega_*}
}
\end{equation}
\end{thm}

\begin{proof}
The canonical representation $(j_{M_L},j_B)$ is the composition of the universal Toeplitz representation $(\psi,\pi)$ of $M_L$ in $\To(M_L)$ with the quotient map $Q$, and in particular $j_B=Q\circ\pi$.  Since $B$ is separable, \cite[Theorem~4.4]{P} says that the homomorphism $\pi:B\to\To(M_L)$ induces an isomorphism $\pi_*:K_i(B)\to K_i(\To(M_L))$, and since $\rho:B\to \ker Q$ is an isomorphism onto a full corner, $\rho_*$ is an isomorphism.
So splicing the commutative diagram of Proposition~\ref{prop-square} into \eqref{eq-les2} gives the result.
\end{proof}

\section{Endomorphisms arising from dilation matrices}

Throughout this section, $d$ is an integer $\geq 2$ and $A\in M_d(\Z)$ is an
integer dilation matrix, by which we mean that all the complex eigenvalues
$\lambda$ of $A$ satisfy $|\lambda|>1$. We consider the surjective endomorphism
$\sigma_A$ of $\T^d$ defined by  $\sigma_A(e^{2\pi ix})=e^{2\pi iAx}$ for
$x\in\R^d$, which has $|\ker \sigma_A|=|\det A|$, and the associated Exel system
$(C(\T^d),\alpha_A, L)$,  
  where $\alpha_A$ is the endomorphism of $C(\T^d)$ given by $\sigma_A$.

We start by showing that $C(\T^d)\rtimes_{\alpha_A,L}\N=\O(M_L)$ is simple and purely infinite. We deduce simplicity from results of Exel and Vershik \cite{EV} on crossed products by endomorphisms, and pure infiniteness from results of Katsura \cite{K4} on the $C^*$-algebras of topological graphs. So we need to note that the map $f\mapsto N^{1/2}f$ is an isomorphism of the bimodule $M_L=C(K)_L$ onto the bimodule of the topological graph $E$ with $E^0=\T^d$, $E^1=\T^d$, $r=\id$ and $s=\sigma_A$, and hence the crossed product $C(\T^d)\rtimes_{\alpha, L}\N=\O(M_L)$ can also be viewed as the $C^*$-algebra $C^*(E)$ studied in \cite{K1, K4}.

We need the following lemma on the operator norms of $A^n$ acting on $\R^d$.

\begin{lemma}\label{lem-dilation}
We have $\|A^{-n}\|\to 0$ as $n\to\infty$.
\end{lemma}

\begin{proof}
For a real matrix $B$, the operator norms of $B$ in $B(\R^d)$ and $B(\C^d)$ coincide (the $C^*$-identities imply that both are equal to the square root of the largest eigenvalue of $B^TB$). So we may as well work over $\C$, and then there exists $P\in GL_d(\C)$ such that $P^{-1}A^{-1}P$ is in Jordan canonical form. Thus $P^{-1}A^{-1}P$ has the form $D+N$ where $D$ is diagonal, $N$ is nilpotent with $N^{d}=0$, and $D$ and $N$ commute. The entries of $D$ are the reciprocals of the eigenvalues of $A$, so 
\[
\|D\|=\max\{|\lambda^{-1}|: \text{ $\lambda$ is an eigenvalue of $A$}\}< 1,
\]
and $\|N\|\leq 1$ because $N$ is a truncated shift. Since $\|A^{-n}\|\leq\|P\|\|P^{-1}\|\|(D+N)^n\|$, it suffices to show that $\|(D+N)^n\|\to 0$ as $n\to\infty$.

Since $D$ and $N$ commute and $N^{d}=0$, for $n\geq d$ the binomial theorem gives
\[
\|(D+N)^n\|=\left\| \sum_{k=0}^{d-1}\binom{n}{k}D^{n-k}N^k\right\|\leq\|D\|^{n-d+1}\sum_{k=0}^{d-1}\binom{n}{k}\|D^{d-1-k}\|\|N\|^k,
\]
and since $\|D^{d-1-k}\|\leq \|D\|^{d-1-k}\leq 1$ for $0\leq k\leq d-1$, we have
\[
\|(D+N)^n\|\leq\|D\|^{-d+1}\|D\|^{n}\sum_{k=0}^{d-1}\binom{n}{k}=\|D\|^{-d+1}\|D\|^{n}f(n)
\]
where $f$ is a polynomial of degree $d-1$.  But $\|D\|^nf(n)=\exp(n\ln\|D\|)f(n)\to 0$ as $n\to\infty$ because $\ln\|D\|<0$, and the lemma follows.
\end{proof}

\begin{prop}\label{prop-simple}
The Cuntz-Pimsner algebra $\O(M_L)$ is simple and purely infinite.
\end{prop}

\begin{proof} 
We show that $\O(M_L)$ is simple using \cite[Theorem~11.2]{EV}, which says that $C(\T^d)\rtimes_{\alpha,L}\N$ is simple if and only if $\sigma_A$ is irreducible.
We recall from  \cite[\S11]{EV} that $x,y\in \T^d$ are
\emph{trajectory-equivalent}, written $x\sim y$,  if
there are $n,m\in\N$ such that $\sigma_A^n(x) = \sigma_A^m(y)$, and a subset $Y\subseteq \T^d$ is \emph{invariant} if $x\sim y \in Y$
implies that $x\in Y$;  $\sigma_A$ is \emph{irreducible} if the only
closed invariant sets are $\emptyset$ and $\T^d$.

Let $Y$ be a non-empty closed invariant subset of $\T^d$, and pick a point $e^{2\pi i y}\in Y$. We need to show that $Y=\T^d$. Fix $e^{2\pi i z}\in \T^d$. Since the unit cube in $\R^d$ has diameter $\sqrt{d}$, \aah{for every $n\in \N$} we can find $k_n\in \Z^d$ such that $|A^nz-(y+k_n)|\leq \sqrt{d}$. Then $x_n:=A^{-n}(y+k_n)$ has $\sigma_A^n(e^{2\pi ix_n})=e^{2\pi iA^nx_n}=e^{2\pi iy}\in Y$, and invariance implies that $e^{2\pi ix_n}\in Y$ also. Lemma~\ref{lem-dilation} implies that 
\[
|z-x_n|\leq\|A^{-n}\||A^nz-(y+k_n)|\leq\|A^{-n}\|\sqrt{d}\to 0\text{\ as $n\to\infty$},
\]
so $x_n\to z$ in $\R^d$ and $e^{2\pi i x_n}\to e^{2\pi i z}$.  Since $Y$ is closed, this implies that $e^{2\pi i z}\in Y$, as required. Thus $\sigma_A$ is irreducible, and $\O(M_L)$ is simple.

To show that $\O(M_L)$ is purely infinite we  realise $\O(M_L)=C(\T^d)\rtimes_{\alpha, L}\N$ as $C^*(E)$ with $E=(\T^d, \T^d,\id, \sigma_A)$. Since $C^*(E)=\O(M_L)$ is simple, $E$ is minimal by \cite[Proposition~1.11]{K4}.
So by \cite[Theorem~A]{K4} it suffices to prove that $E$ is contracting at some vertex $v_0\in E^0$ in the sense of Definition~2.3 of \cite{K4}; we will show that $E$ is contracting at $v=(1,1,\dots,1)$.   First, we need to see that the positive  orbit $\{z:\sigma_A^n(z)=v\}$ of $v$ is dense in $E^0=\T^d$.   The positive orbit of $v$ contains all points of the form $e^{2\pi i A^{-n}k}$ for $n\in N$ and $k\in \Z^d$, and it follows from our proof of the irreducibility of $\sigma_A$ above (with $y=0$) that this positive orbit is dense in $E^0$.  

Second, we fix a neighbourhood $V$ of $v$; we need to show that $V$ contains a contracting open set $W$ (see \cite[Definition~2.3]{K4}). For this, it suffices to find a open neighbourhood $W$ of $v$ such that $W\subset V$ and $\overline{W}\subsetneq \sigma_A^k(W)$ for some $k\geq 1$. By Lemma~\ref{lem-dilation} we can choose $k$ such that $\|A^{-k}\|<1$. Then for every $\epsilon>0$ and every $x$ in the closed unit ball $\overline{B(0,\epsilon)}$ in $\R^d$, we have $|A^{-k}x|<\epsilon$, so $x=A^k(A^{-k}x)$ belongs to $A^k(B(0,\epsilon))$. Thus $\overline{B(0,\epsilon)}\subset A^k(B(0,\epsilon))$. The inequality $\|A^kA^{-k}\|\leq \|A^k\|\,\|A^{-k}\|$ implies that $\|A^k\|>1$, so for every $\epsilon>0$ there exists $y\in B(0,\epsilon)$ such that $|A^ky|>\epsilon$, and $\overline{B(0,\epsilon)}\subsetneq A^k(B(0,\epsilon))$. If $\epsilon$ is small enough to ensure that $x\mapsto e^{2\pi ix}$ is one-to-one on $A^k(B(0,\epsilon))$, then $W:=\{e^{2\pi i x}:x\in B(0,\epsilon)\}$ satisfies $\overline{W}\subsetneq \sigma_A^k(W)$, and by taking $\epsilon$ smaller still we can ensure that $W\subset V$. Thus $E$ is contracting, and the result follows from \cite[Theorem~A]{K4}. 
\end{proof}

We now want to calculate the $K$-theory of $C(\T^d)\rtimes_{\alpha_A,L}\N=\O(M_L)$, and we aim to use Theorem~\ref{thm-stes}. To do this, we need descriptions of $K_*(C(\T^d))$ and the map $\Omega_*$. 

\begin{lemma}\label{lem-omega*}
Suppose that $(B,\alpha,L)$ is an Exel system  with $B$ commutative, that $M_L$ admits an orthonormal basis $\{m_j:0\leq j<N\}$, and that $\Omega:B\to M_N(B)$ is the homomorphism described in Lemma~\ref{lem-Omega}. Then $\Omega_*\circ\alpha_*$ is multiplication by $N$ on both $K_0(B)$ and $K_1(B)$.
\end{lemma}

\begin{proof}
We know from Lemma~\ref{lem-Omega}(3) that $\Omega\circ \alpha(a)=a1_n$. If $b\in M_r(B)$, then $(\Omega\circ\alpha)_r(b)$ is the $N\times N$ block matrix which has $0$s off the diagonal and $\alpha_r(b)$ down the diagonal. If we view $(\Omega\circ\alpha)_r(b)$ as an element of $M_r(M_N(B))$, as in Lemma~\ref{tocome}, it becomes $b\oplus b\oplus\cdots \oplus b$. Whether $b$ is a projection or a unitary, $[b\oplus\cdots \oplus b]=N[b]$. Thus by Lemma~\ref{tocome}, we have
\[
\Omega_*\circ\alpha_*([b])=(\Omega\circ\alpha)_*([b)]=[(\Omega\circ\alpha)_r(b)]=[b\oplus\cdots \oplus b]=N[b].\qedhere
\]
\end{proof}

Ji proved in \cite{J} that the Chern character is a $\Z/2$-graded ring isomorphism of $K_*(C(\T^d))=K^*(\T^d)$ onto the integral cohomology ring 
\[
\textstyle{H^*(\T^d,\Z):=\bigoplus_{k\in\Z}^\infty H^k(\T^d,\Z)=\bigoplus_{k=0}^d H^k(\T^d,\Z),}
\] 
which in turn is isomorphic as a $\Z$-graded ring to the exterior algebra $\bigwedge^*\Z^d$. Thus the ring $H^*(\T^d,\Z)$ is generated by $H^1(\T^d,\Z)$, which is isomorphic to the set of homotopy classes of continuous functions from $\T^d$ to $\T$, and is the free abelian group generated by the coordinate functions $u_k:z=(z_1,\cdots, z_n)\mapsto z_k$. Since the homomorphism $\alpha_*$ is induced by a continuous map $\sigma_A:\T^d\to \T^d$, the corresponding ring homomorphism on $H^*(\T^d,\Z)$ is the map $\sigma_A^*$, which respects the $\Z$-grading. Thus we can compute $\alpha_*$ on $\bigwedge^*\Z^d$ by working out what $\sigma_A^*$ does on $H^1(\T^d,\Z)$ using the basis $\{e_k:=[u_k]:1\leq k\leq d\}$, and then taking exterior powers. Once we know what $\alpha_*$ is, we can use the formula for $\Omega_*\circ\alpha_*$ in Lemma~\ref{lem-omega*} to work out what $\Omega_*$ is.

\begin{lemma}\label{lem-alpha*1}
With respect to the basis $\{[u_k]\}$, $\alpha_*:\lsp\{[u_k]\}\to \lsp\{[u_k]\}$ is multiplication by the transpose $A^T$ of $A$.
\end{lemma}

\begin{proof}
We have  $\alpha_*([u_k])=[\alpha(u_k)]=[u_k\circ \sigma_A]$. Since
\begin{align*}
u_k\circ\sigma_c(e^{2\pi ix})&=u_k(e^{2\pi iAx})
=e^{2\pi i\sum_j a_{k,j}x_j}=\prod_j e^{2\pi ia_{k,j}x_j}\\
&=\prod_j (e^{2\pi ix_j})^{a_{k,j}}=\prod_ju_j(e^{2\pi ix})^{a_{k,j}},
\end{align*}
we have $u_k\circ\sigma_A=\prod_j u_j^{a_{k,j}}$.  Hence $[u_k\circ \sigma_A]=\sum_j a_{k,j}[u_j]$.
\end{proof}

Since the $0$-graded component is isomorphic to $H^0(\T^d,\Z)$, the free abelian group generated by the connected components, the action of $\alpha_*$ on the $0$-component $\bigwedge^0(\Z)=\Z$ is the identity map. For $n=1$, Lemma~\ref{lem-alpha*1} implies that $\alpha_*=A^T$. For $n> 1$, we use the basis
\[
\Ee_n=\big\{e_J=e_{j_1}\wedge\dots\wedge e_{j_n}:J\subset\{1,\dots,d\}, |J|=n,J=\{j_1<j_2<\dots <j_n\}\big\}
\]
for $\bigwedge^n\Z^d$. For $e_K\in \Ee_n$, we write $K'=\{1,\dots, d\}\setminus K$. With $K$ and $K'$ listed in increasing order as $K=\{k_1<\dots<k_n\}$ and $K'=\{k_{n+1}<\dots<k_d\}$, we let $\tau_K$ be the permutation $i\mapsto k_i$ for $1\leq k\leq d$. For subsets $K,J$ of the same size, we write $A_{K,J}$ for the submatrix of $A$ whose entries belong to the rows in $K$ and the columns in $J$. The following Lemma is essentially Lemma~1 of \cite[Chapter~5]{N}; we have included a short proof because the conventions of \cite{N} are different (matrices act on the right of vector spaces, for example).

\begin{lemma}
Let $1\leq n\leq d$. The matrix $C_n$ of $\alpha_*|:=\bigwedge^n \Z^d\to \bigwedge^n \Z^d$  with respect to the basis $\Ee_n$ has $(J,K)$ entry $\det A_{K,J}$.
\end{lemma}

\begin{proof}
Fix $e_K\in\Ee_n$ with $K=\{k_1<\dots <k_n\}$.  Then
\begin{align*}
(\bigwedge{}^n A^T)(e_K)&=(\bigwedge{}^n A^T)(e_{k_1}\wedge\dots\wedge e_{k_n})\\
&=A^Te_{k_1}\wedge\dots\wedge A^T e_{k_n}\\
&=\sum_{m_1=1,\dots, m_n=1}^d a_{k_1,m_1}\dots a_{k_n, m_n}(e_{m_1}\wedge\dots\wedge e_{m_n})\\
&=
\sum_{e_J\in\Ee_n}\sum_{\{m_1,\dots,m_n\}=J}a_{k_1,m_1}\dots a_{k_n, m_n}(e_{m_1}\wedge\dots\wedge e_{m_n})\\
&=\sum_{e_J\in\Ee_n}\sum_{\sigma\in S_n} a_{k_1,\sigma(j_1)}\dots a_{k_n, \sigma(j_n)}(e_{\sigma(j_1)}\wedge\dots\wedge e_{\sigma(j_n)})\\
&=\sum_{e_J\in\Ee_n}\sum_{\sigma\in S_n}(-1)^{\text{deg}\sigma} a_{k_1,\sigma(j_1)}\dots a_{k_n, \sigma(j_n)}(e_{j_1}\wedge\dots\wedge e_{j_n})\\
&=\sum_{e_J\in\Ee_n}(\det A_{K,J})e_J.\qedhere
\end{align*}
\end{proof}

We are now ready to compute the matrix $B_n$ of $\Omega_*$ on $\bigwedge\Z^d$ with respect to the same basis \aah{$\Ee_n$}. The answer must, of course, be an integer matrix. But Lemma~\ref{lem-omega*} implies that $C_n$ is invertible as a real matrix, and hence if we can find matrices $B_n$ such that $B_nC_n=N1_n$, then uniqueness of the real inverse tells us that $B_n$ is the matrix of $\Omega_*$.

\begin{prop}\label{lem-helper2} Let $B_0=|\det A|$, $B_d=\sign(\det A)$, and
\[
B_n=\begin{cases}
\Big( (-1)^{\deg(\tau_K\tau_L)} \det(A_{K',L'})\Big)_{K,L}&\text{if $\det A>1$;}\\
-\Big( (-1)^{\deg(\tau_K\tau_L)} \det(A_{K',L'})\Big)_{K,L}&\text{if $\det A<-1$.}
\end{cases}
\]

\smallskip
\textnormal{(1)} Then $B_nC_n=|\det A| 1$ where $1$ is the $\binom{d}{n}\times \binom{d}{n}$ identity matrix.

\smallskip
\textnormal{(2)} We have $1-B_0= 1-|\det A|<0$, $\det (1- B_n)\neq 0$ for $1\leq n<d$, and 
\[1-B_d=\begin{cases}
0&\text{if $\det A>1$}\\
2&\text{if $\det A<-1$.}
\end{cases}\] 
\end{prop}

\aah{For the proof of Proposition~\ref{lem-helper2} we need the following lemma; its first part appears as equation (5.3.7) in \cite{N}, for example.}

\begin{lemma}\label{lem-helper1}
Fix $n$ satisfying $1\leq n\leq d-1$. 

\smallskip
\textnormal{(1)} If $e_J\in \Ee_n$, then
 \[
 \det A=\sum_{e_K\in\Ee_n} (-1)^{\deg(\tau_K\tau_J)}\det(A_{K,J})\det(A_{K',J'}).
\]

\textnormal{(2)}  If $e_J,e_L\in \Ee_n$ and $L\neq J$, then
\[
\sum_{e_K\in\Ee_n}(-1)^{\deg(\tau_K\tau_J)}\det(A_{K,J})\det(A_{K',L'})=0.
\]
\end{lemma}

\begin{proof}
\noindent (1) Fix $e_J\in\Ee_n$.
We have
\begin{align}
\det A
&=\sum_{\sigma\in S_d}(-1)^{\deg\sigma}a_{\sigma(1),1}\dots a_{\sigma(d),d}\notag\\
&=(-1)^{\deg\tau_J}\sum_{\sigma\in S_d}(-1)^{\deg\sigma}a_{\sigma(1),j_1}\dots a_{\sigma(d),j_d}\notag\\
\intertext{which, by reordering the sum according to the image of $I_n:=\{1,\dots,n\}$ under $\sigma$, is }
&=(-1)^{\deg\tau_J}\sum_{e_K\in\Ee_n}\sum_{\{\sigma:\sigma(I_n)=K\}}(-1)^{\deg\sigma}a_{\sigma(1),j_1}\dots a_{\sigma(d),j_d}.\label{eq-headache}
\end{align}
Note that for fixed $\sigma\in S_n$ such that $\sigma(I_n)=K$ we have
 \[
 \sigma=(\sigma_K\times\sigma_{K'})\circ \tau_K
 \]
 where $\sigma_K(k_i):=\sigma(i)$ and $\sigma_{K'}(k_l):=\sigma(l)$. So
\begin{align*}
\eqref{eq-headache}
&=
(-1)^{\deg\tau_J}\sum_{e_K\in\Ee_n}\sum_{\{\sigma:\sigma(I_n)=K\}}(-1)^{\deg\tau_K}(-1)^
{\deg(\sigma_K\times\sigma_{K'})}a_{\sigma_K(k_1),j_1}\dots a_{\sigma_K(k_n),j_n}\cdot\\
&\hskip7cm\cdot a_{\sigma_{K'}(k_{n+1}),j_{n+1}}\dots a_{\sigma_{K'}(k_{d}),j_{d}}\\
&=
(-1)^{\deg\tau_J}\sum_{e_K\in\Ee_n}(-1)^{\deg\tau_K}\sum_{\alpha\in S_K,\beta\in S_{K'}}
(-1)^{\deg\alpha} a_{\alpha(k_1),j_1}\dots a_{\alpha(k_n),j_n}\cdot\\
&\hskip7cm \cdot (-1)^{\deg\beta} a_{\beta(k_{n+1}),j_{n+1}}\dots a_{\beta(k_{d}),j_{d}}\\
&=
\sum_{e_K\in\Ee_n}(-1)^{\deg(\tau_K\tau_J)}\det(A_{K,J})\det(A_{K',J'}).
\end{align*}

(2) If $L\neq J$ then $L'\neq J'$ and $L'\cap J\neq \emptyset$. Consider the matrix $D$ whose entries are those of $A$ except that the $L\setminus J$ columns of $D$ have been replaced by copies of the $J\setminus L$ columns of $A$.  Thus $\det D=0$.  Note that $A_{K,J}$ and $D_{K,L}$ have the same columns up to permutation, so $\det(A_{K,J})=\pm\det(D_{K,L})$.  For every $K$ we have $D_{K',L'}=A_{K',L'}$, so using (1) we get
\begin{align*}
\sum_{e_K\in\Ee_n}&(-1)^{\deg(\tau_K\tau_J)}\det(A_{K,J})\det(A_{K',L'})\\
&=\pm\sum_{e_K\in\Ee_n}(-1)^{\deg(\tau_K\tau_J)}\det(D_{K,L})\det(D_{K',L'})=\det D=0.\qedhere
\end{align*}
\end{proof}

\begin{remark}
In \cite[page~92]{N}, it is observed that the coefficient $(-1)^{\deg(\tau_K\tau_J)}$ can be realised as the product $\prod_{i=1}^n(-1)^{j_i+k_i}$. To see this, first observe that $(-1)^{\deg(\tau_J)}=\prod_{i=1}^n(-1)^{j_i-i}$ (because $j_n-n$, for example, is the number of transpositions required to move $j_n$ to its correct place in $J'$ without changing the ordering of $J'$), and then $(-1)^{\deg(\tau_K\tau_J)}=\prod_{i=1}^n(-1)^{(j_i-i)+(k_i-i)}$.
\end{remark}

\begin{proof}[Proof of Proposition~\ref{lem-helper2}]
Say $\det A>1$.  Then the $(J,L)$ entry of $C_n B_n$ is
\[
\sum_{e_K\in\Ee_n}\det (A_{K,J})(-1)^{\deg(\tau_K\tau_L)}\det(A_{K',L'})
\]
which, by Lemma~\ref{lem-helper1}, equals $\delta_{J,L}(\det A)1$.   If $\det A<-1$ the same calculation gives $-\delta_{J,L}(\det A)1=\delta_{J,L}|\det A|1$. Thus $C_n B_n=|\det A|1=B_nC_n$. This gives (1).

\smallskip

The statements in (2) about $B_0$ and $B_d$ are immediate, so we suppose $1\leq n\leq d-1$. To compute $\det(I-B_n)$ we work over $\C$, and choose a basis for $\C^d$ such that $A$ is upper-triangular. We claim that if $J=\{j_1<\dots<j_n\}>K=\{k_1<\dots<k_n\}$ in the lexicographical order, then $\det (A_{J,K})=0$. If $J>K$ then there exists $m$ such that $j_i=k_i$ for $i<m$ and $j_m>k_m$.  Since $A$ is upper-triangular, $j_m>k_m$ implies $a_{j_m,k_m}=0$.  Moreover, $j_{n-m+1}>\dots>j_{m+1}>j_m>k_m$,  so $A_{J,K}$ has the form
\[A_{J,K}=\begin{pmatrix}
U&*\\
0&V
\end{pmatrix}
\] where $U$ is an $(m-1)\times (m-1)$ upper-triangular matrix, and $V$ is a square matrix with  the first column consisting of zeros.  Thus $\det (A_{J,K})=0$, as claimed.  So if we order $\Ee_n$ with the lexicographic order, then the matrix  $(\det(A_{K,J}))_{J,K}$ of $\alpha_*|=\bigwedge^n A^T$ is lower-triangular.  Hence its inverse $(\det A)^{-1}B_n$ is also lower-triangular, and so is $B_n$.  The diagonal entries of $B_n$ are $\det(A_{K',K'})=\prod_{k\in K'}a_{k,k}$; since each $a_{k,k}$ is an eigenvalue of $A$, we have $|a_{k,k}|>1$, and each diagonal entry of $B_n$ has absolute value greater than $1$. Since $B_n$ is lower-triangular, it follows that $\det(1- B_n)\neq 0$.
\end{proof}

\begin{thm}\label{thm-example}
Let $A$ be a dilation matrix $A\in GL_d(\Z)$ with $d\geq 1$, and define $B_n$ as in Proposition~\ref{lem-helper2}. Let $M_L$ be the bimodule for the Exel system $(C(\T^d),\alpha_A,L)$ and for which $C(\T^d)\rtimes_{\alpha_A,L}\N=\O(M_L)$.

\smallskip
\textnormal{(1)} If $\det A>1$ and $d$ is odd, then
\begin{align*}
 K_0(\O(M_L))&=
 {\textstyle\big(\bigoplus_{\textnormal{$n$ even, $n<d$}}
 \coker(1-B_n)\big)\oplus\Z},\ \text{and}\\
 K_1(\O(M_L))&=
 {\textstyle \bigoplus_{\textnormal{$n$ odd, $n\leq d$}}
 \coker(1-B_n).}
 \end{align*}
 If $\det A>1$ and $d$ is even, then
 \begin{align*}
  K_0(\O(M_L))&=
 {\textstyle\bigoplus_{\textnormal{$n$ even, $n\leq d$}}
 \coker(1-B_n)},\ \text{and}\\
 K_1(\O(M_L))&=
 \big({\textstyle \bigoplus_{\textnormal{$n$ odd, $n< d$}}
 \coker(1-B_n)\big)\oplus\Z.}
 \end{align*}

\textnormal{(2)} If $\det A<-1$, then
 \begin{align*}
  K_0(\O(M_L))&=
 {\textstyle\bigoplus_{\textnormal{$n$ even, $n\leq d$}}
 \coker(1-B_n)},\ \text{and}\\
 K_1(\O(M_L))&=
{\textstyle \bigoplus_{\textnormal{$n$ odd, $n\leq d$}}
 \coker(1-B_n).}
 \end{align*}
 
\end{thm}

\begin{proof}  We identify
\[\textstyle{K_1(C(\T^d))\cong\bigoplus_{\textnormal{$n$ odd, $n\leq d$}}\bigwedge^n\Z^d\quad\text{and}\quad K_0(C(\T^d))\cong\bigoplus_{\textnormal{$n$ even, $n\leq d$}}\bigwedge^n\Z^d.}\]

\smallskip
\aah{Suppose that $\det A>1$. By Lemma~\ref{lem-omega*}, $(\Omega\circ\alpha)_*$ is multiplication by  $|\det A|$, and by Proposition~\ref{lem-helper2}(1) the matrix $C_n$ of $\alpha_*|$ has inverse $|\det A|^{-1}B_n$; it follows that the map} $\id-\Omega_*$ appearing in Diagram~\ref{eq-les} is \[\textstyle{\bigoplus_{\textnormal{even $n\leq d$}}(1-B_n)\quad\text{and}\quad \bigoplus_{\textnormal{odd\ } n\leq d}(1-B_n)}\] on  $K_0(C(\T^d))$ and $K_1(C(\T^d))$, respectively. By Proposition~\ref{lem-helper2}(2), each $1-B_n$ with $n<d$ is injective, and $1-B_d=0$.

Suppose that $d$ is odd. Then $\bigoplus_{\textnormal{even $n\leq d$}}(1-B_n)$ is injective and  \[\ker \big(\textstyle{\bigoplus_{\textnormal{odd $n\leq d$}}}(1-B_n)\big)=\ker(1-B_d)=\Z.\] Thus Diagram~\ref{eq-les} gives
\[K_1(\O(M_L))\cong
 \textstyle{\bigoplus_{\textnormal{$n$ odd, $n\leq d$}}}
 \coker(1-B_n)\]
 and an exact sequence
\begin{equation*} 
\xymatrix{
0\ar[r]&\textstyle{\bigoplus_{\textnormal{$n$ even, $n< d$}}}\coker(1-B_n)
\ar[r]&K_0(\O(M_L))\ar[r]&\Z\ar[r]&0.
}
\end{equation*}
Since $\Z$ is free this sequence splits, and the formula for $K_0$ follows. 

The proof for even $d$ is similar.

For part (2), we just note that Proposition~\ref{lem-helper2}(2) implies that   
\[{\textstyle\bigoplus_{\textnormal{$n$ even, $n\leq d$}}}(1-B_n)\quad\text{and}\quad {\textstyle\bigoplus_{\textnormal{$n$ odd, $n\leq d$}}}(1-B_n)\] are injective, and the result follows.
\end{proof}

For small $d$ we can identify the $B_n$ in more familiar terms. Both $B_0$ and $B_d$ are just numbers (or rather, multiplication by those numbers on $\Z$). Next we have:

\begin{prop}\label{B1Bd-1}
For every $d$ we have $B_1=|\det A|(A^T)^{-1}$. If we list the basis for $\bigwedge^{d-1}\Z^d$ as $f_k:=e_{\{1,\cdots,d\}\setminus\{k\}}$, then $B_{d-1}$ is the matrix with $(k,l)$ entry $(-1)^{k+l}a_{k,l}$ (if $\det A>0$) or $(-1)^{k+l+1}a_{k,l}$ (if $\det A<0$).
\end{prop}

\begin{proof}
For each singleton set $\{k\}$, the permutation $\tau_{\{k\}}$ is the cycle which pulls $k$ to the front and moves the elements $1,\cdots k-1$ to the right, which has degree $k$. The complements $\{k\}'$ are the sets $\hat k:=\{1,\cdots,k\}\setminus\{k\}$, and the number
\[
(-1)^{\deg(\tau_K\tau_L)} \det A_{\{k\}',\{l\}'}=(-1)^{\deg \tau_K+\deg\tau_L}\det A_{\hat k,\hat l}=(-1)^{k+l}\det A_{\hat k,\hat l}
\]
is the $(l,k)$ entry in $(\det A)A^{-1}$, and the $(k,l)$ entry in $(\det A)(A^T)^{-1}$. The extra minus sign in the formula for $B_1$ when $\det A<0$ shows that $B_1$ is $(|\det A|)(A^T)^{-1}$.

The $(k,l)$ entry in the matrix of $B_{d-1}$ with respect to the basis $\{f_k\}$ is the $(\hat k,\hat l)$ entry in the matrix with respect to the basis $\Ee_{d-1}$. For $K=\hat k$, $\tau_K$ is the cycle which moves $k$ to the back and the last $d-k$ terms one forward, which has degree $d-k+1$. Since $A_{(\hat k)',(\hat l)'}$ is the $1\times 1$ matrix with entry $a_{k,l}$, we have
\[
(-1)^{\deg(\tau_K\tau_L)} \det A_{(\hat k)',(\hat l)'}=(-1)^{(d-k+1)+(d-l+1)}a_{k,l}=(-1)^{2(d+1)-(k+l)}a_{k,l}=(-1)^{k+l}a_{k,l}.
\]
This immediately gives the result for $\det A>0$, and for $\det A<0$, the extra minus sign in the formula for $B_{d-1}$ means we need to replace $(1)^{k+l}$ by $(-1)^{k+l+1}$. 
\end{proof}

We can now sum up our results for small $d$: Corollary~\ref{cord=1} is well-known, as we observed in the introduction, but  Corollary~\ref{calcd=2} was a bit of a surprise. 

\begin{cor}\label{cord=1}
Suppose $N$ is a non-zero integer, and consider the Exel system $(C(\T),\alpha_N,L)$ associated to the covering map $z\mapsto z^N$.

\smallskip
\textnormal{(1)} If $N>1$, then $K_0(C(\T)\rtimes_{\alpha_N,L}\N)=(\Z/(N-1)\Z)\oplus\Z$ and $K_1(C(\T)\rtimes_{\alpha_N,L}\N)=\Z$.

\smallskip
\textnormal{(2)} If $N<-1$, then $K_0(C(\T)\rtimes_{\alpha_N,L}\N)=\Z/(N-1)\Z$ and $K_1(C(\T)\rtimes_{\alpha_N,L}\N)=\Z/2\Z$.
\end{cor}

\begin{cor}\label{calcd=2}
Suppose that $A=(a_{ij})\in M_2(\Z)$ is a dilation matrix. Then
\[
K_0(C(\T^2)\rtimes_{\alpha_A,L}\N)=
\begin{cases}
\Z/(|\det A|-1)\Z\oplus\Z&\text{if $\det A>1$}\\
(\Z/(|\det A|-1)\Z)\oplus (\Z/2\Z)&\text{if $\det A<-1$,}
\end{cases}
\]
and
\[
K_1(C(\T^2)\rtimes_{\alpha_A,L}\N)=
\begin{cases}
\Z\oplus\coker\bigg(\begin{matrix}1-a_{11}&a_{12}\\a_{21}&1-a_{22}\end{matrix}\bigg)&\text{if $\det A>1$}\\
\coker\bigg(\begin{matrix}1+a_{11}&-a_{12}\\-a_{21}&1+a_{22}\end{matrix}\bigg)&\text{if $\det A<-1$.}
\end{cases}
\] 
\end{cor}

\begin{proof}
The statement about $K_0$ follows immediately from Theorem~\ref{thm-example}. For $K_1$, we use the description of $B_1=B_{2-1}$ in Proposition~\ref{B1Bd-1}. (If we had used the description of $B_1$ as $|\det A|(A^T)^{-1}$, we would have got a different matrix, because we would then be calculating it with respect to the basis $\{e_1,e_2\}$ rather than $\{f_1,f_2\}=\{e_2,e_1\}$. However, the two matrices are conjugate in $M_2(\Z)$, and hence have isomorphic cokernels.)
\end{proof}

We now look at the implications of these results for some concrete examples of  dilation matrices. The first two were used in \cite{P2} to provide examples of projective multi-resolution analyses.

\begin{examples} (1) The matrix $A=\big(\begin{smallmatrix}0&1\\2&0\end{smallmatrix}\big)$ has $\det A=-2<-1$. So $K_1(C(\T^2)\rtimes_{\alpha_A,L}\N)$ is the cokernel of \aah{$\big(\begin{smallmatrix}1&-1\\-2&1\end{smallmatrix}\big)$}; since this matrix has determinant $-1$, it is invertible over $\Z$, and we have
\[
K_0(C(\T^2)\rtimes_{\alpha_A,L}\N)=\Z/2\Z \text{ and } K_1(C(\T^2)\rtimes_{\alpha_A,L}\N)=0.
\]
These $K$-groups are the same as those of $\O_3$, but the class of the identity is different. To see the last statement, note that the class $[1]$ of the identity in $K_0(C(\T^2))$ is the image of $1\in \Z=\bigwedge^0\Z^2$, and when $|\det A|=2$, $1-B_0$ is invertible, so $[1]$ belongs to the range of $\id-\Omega_*$. Thus the class of the identity $1_{C(\T^2)\rtimes\N}=j_{C(\T)}(1)$ in $K_0(C(\T^2)\rtimes_{\alpha_A,L}\N)$ is $0$. For $\O_3$, on the other hand, $[1]$ is the generator of $K_0(\O_3)$.

\smallskip
(2) The matrix $A=\big(\begin{smallmatrix}1&1\\-1&1\end{smallmatrix}\big)$ has $\det A=2>1$. So Corollary~\ref{calcd=2} implies that 
\[
K_0(C(\T^2)\rtimes_{\alpha_A,L}\N)=\Z \text{ and } K_1(C(\T^2)\rtimes_{\alpha_A,L}\N)=\Z.
\]

\smallskip
(3) The matrix $A=\big(\begin{smallmatrix}2&1\\-1&2\end{smallmatrix}\big)$ has $\det A=5>1$. Thus 
\[
K_0(\O(M_L))=\Z/4\Z\oplus\Z \text{ and } K_1(C(\T^2)\rtimes_{\alpha_A,L}\N)=\Z\oplus (\Z/2\Z).
\]

\smallskip
(4) The matrix $A=\big(\begin{smallmatrix}2&-1\\1&-3\end{smallmatrix}\big)$ has determinant $-5$, and
\[
K_0(C(\T^2)\rtimes_{\alpha_A,L}\N)=(\Z/4\Z)\oplus (\Z/2\Z)\text{ and } K_1(C(\T^2)\rtimes_{\alpha_A,L}\N)=\Z/5\Z.
\]
\end{examples}

No, we don't see any obvious pattern either.

\end{document}